\documentclass[11pt]{article}

\usepackage{epsf,rotating,amssymb}

\def\al{\alpha}
\def\be{\beta}
\def\la{\lambda}
\def\om{\omega}
\def\ip{\lrcorner\,}
\def\ds{\oplus}
\def\map{\rightarrow}
\def\bq{\begin{equation}}
\def\eq{\end{equation}}
\def\li{\langle}
\def\ri{\rangle}
\def\rc{{\mathbb R}}

\def\fg{{\mathfrak g}}
\def\fh{{\mathfrak h}}
\def\ff{{\mathfrak f}}
\def\ti{\tilde}

\def\Si{\Sigma}
\def\Ga{\Gamma}
\def\La{\Lambda}
\def\ss{\subset}

\begin{document}

\begin{center}
\vskip 2cm {\LARGE Some title containing the words ``homotopy'' and
``symplectic'', e.g. this one\footnote{Based on a talk given at
Poisson 2000, CIRM, Marseille, June 2000.\\
The research for this paper was supported by the European Postdoctoral
Institute (EPDI).\\email: \texttt{severa@sophia.dtp.fmph.uniba.sk}}}\\
\vskip 15mm
{\Large Pavol \v Severa}\\
\vskip 3mm
{\it Dept. of Theor. Physics\\
Faculty of Mathematics and Physics\\
Comenius University\\
842 15 Bratislava, Slovakia}

\end{center}
\vskip 15mm

There is an infinite series of notions, starting with symplectic
manifolds ($n=0$), Poisson manifolds ($n=1$) and Courant algebroids
($n=2)$; there seems to be no name for higher $n$'s, so let us call
the $n$'th term $\Si_n$-manifolds. Their overview is in Table 1 at
the end of the paper. Except for the non-standard terminology, this
table is well known (the connection with variational problems may
be an exception); nevertheless, it seems interesting to write a
short informal review. We'll be mostly concerned with homotopy (or
integration) of $\Si_n$-manifolds. The only non-trivial column is
the one about quantization; for this reason, it won't be mentioned
anymore.

The paper is based on a straightforward use of the basic idea of
Sullivan's Rational homotopy theory \cite{inf} in differential geometry.
Its connection with symplectic geometry is from \cite{aksz}.

\section{Integration of Lie algebroids (after Dennis Sullivan)}

Let us begin with a simple construction of a groupoid $\Ga$ out of a Lie
algebroid $A\map M$. Intuitively, $A$ consists of infinitesimal
morphisms of $\Ga$; to get all the morphisms, we have to compose them
along curves. Thus, consider a Lie algebroid morphism $TI\map A$ ($I$ is
an interval), covering a curve $I\map M$.
 $$\epsfxsize 7cm
\epsfbox{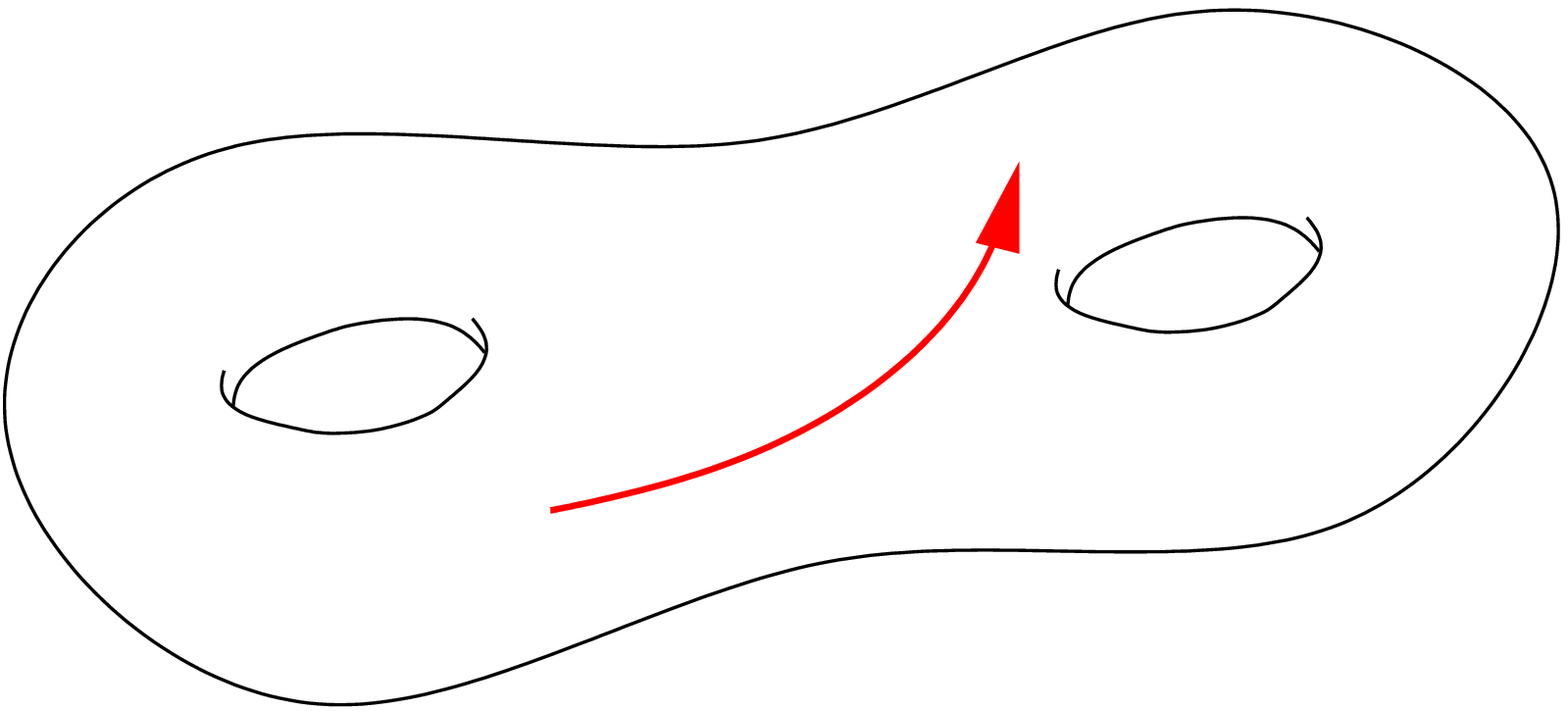}$$
 It gives us a morphism in $\Ga$ between the
endpoints of the curve. Two such morphisms $TI\map A$, with the same
endpoints in $M$, give the same morphism in $\Ga$, if they can be
connected by a morphism $TD\map A$ ($D$ is a disk, or better, to avoid
problems with smoothness, a square with two opposite sides shrunk to
points):
 $$\epsfxsize 7cm \epsfbox{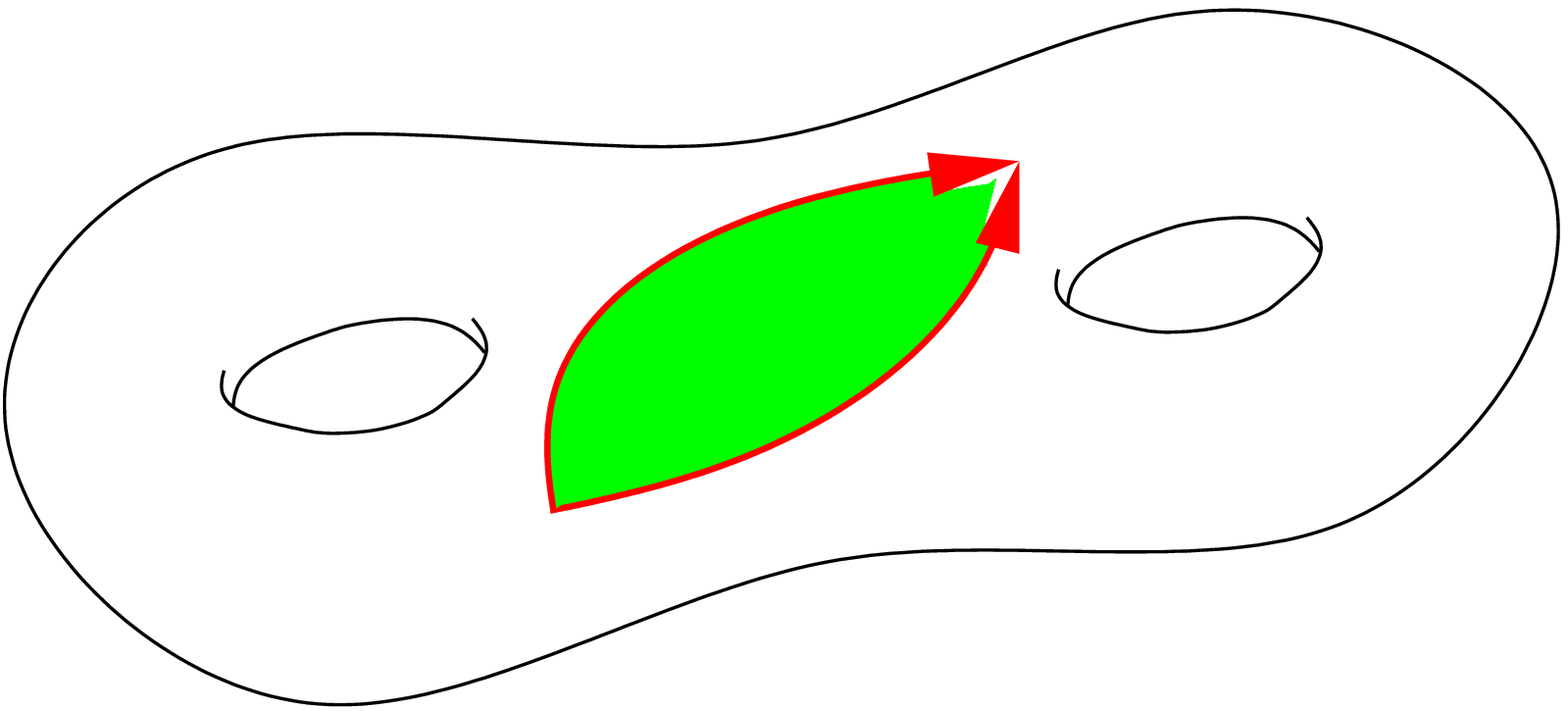}$$
 The composition in
$\Ga$ is just concatenation (again, there is a problem with smoothness
at the joint of the two paths---see below for a proper definition).

Let us notice how similar this construction of $\Ga$ is to the
definition of fundamental groupoid. In what follows, we'll be looking at
generalizations of Lie algebroids with non-trivial higher homotopies,
and at their connections with symplectic geometry.

Finally, a little warning: the $\Ga$ we have defined need not be a
manifold. The problem is that very large morphisms $TD\map A$ may
connect very close morphisms $TI\map A$; thus, the equivalence classes
may not be closed. On the other hand, when we only consider small
morphisms, this construction gives a nice local Lie groupoid. (In the
case when $\Ga$ happens to be a manifold, it is a Lie groupoid with
1-connected fibres.)
 Here is the moral: {\em beware of large disks},
or, more generally, {\em beware of large homotopies}.

For completeness, let us give a definition of $\Ga$ free of the
smoothness problems. Two Lie algebroid morphisms $f_{1,2}:TI\map A$ will
be called {\em homotopic rel boundary (r.b.)} if there is a morphism
$g:T\square\map A$ (where $\square$ is a square) such that the
restrictions of $g$ to $TJ_{1,2}$, where $J_{1,2}$ are the horizontal
sides of $\square$, are constant maps, while the restrictions to
$TI_{1,2}$, where $I_{1,2}$ are the vertical sides, are the morphisms
$f_{1,2}$. $\Ga$ is the space of homotopy classes r.b. of morphisms
$TI\map A$. A morphism $h_3:TI\map A$ will be a composition of two
morphisms $h_{1,2}:TI\map A$ if there is a morphism $k:T\triangle\map A$
(where $\triangle$ is a triangle) s.t. the restriction of $k$ to the
sides of $\triangle$ gives the three $h$'s (with the obvious orientation
of the sides). Any two compositions are homotopic r.b., so the
composition in $\Ga$ is well defined. Any two morphisms $TI\map A$ with a common endpoint in $M$
can be composed, since there is a smooth retraction of $\triangle$ to a pair of its sides.

A closely related and equivalent construction (specialized to the
case of Poisson manifolds) appeared in a paper of Cattaneo and
Felder \cite{cafe}, who also considered the problem of smoothness
of $\Ga$ by describing it as the leaf space of a foliation of an
infinite-dimensional Banach manifold (of $C^1$-morphisms $TI\map
A$). Very recently the problem of smoothness of $\Gamma$ has been
completely solved by Crainic and Fernandes \cite{crfe}.

\section{Some definitions from rational homotopy theory}

Let us follow Sullivan further and make the following definitions: An
{\em N-manifold} (shorthand for ``non-negatively graded supermanifold'')
is a supermanifold with action of the
 multiplicative semigroup
$(\rc,\times)$ such that $-1$ acts as the
 parity operator (i.e. it just
changes the sign of the odd coordinates). When we speak about the degree
of a function (or of a vector field, etc.) on an N-manifold, we mean the
weight of this action, i.e. $f(\la\cdot x)=\la^{\deg f} f(x)$,
$\la\in\rc$. Notice that the degrees of functions are always
non-negative integers and that any function can be approximated (in the
appropriate sense) by finite sums of homogeneous functions (since the
degrees are just exponents of Taylor expansions around $0$ in the
semigroup $\rc$). Finally, an {\em NQ-manifold} is an N-manifold with a
degree-one vector field $Q$ of square 0.

A basic example of an NQ-manifold is $T[1]M$, where $M$ is an ordinary
manifold (the functions on $T[1]M$ are the differential forms on $M$
with their usual degrees, and $Q$ is the de Rham differential). This
will be our basic principle:

{\em Principle: we'll treat NQ-manifolds as if they were of the form
$T[1]M$.}\footnote{This principle is reasonable,
as $T[1]$ is a full and faithful functor from manifolds to NQ-manifolds. It is also
useful to notice that as a functor from N-manifolds to NQ-manifolds, $T[1]$ is the right
adjoint of the corresponding forgetful functor.}

In particular, if $N$ is a manifold and $X$ an NQ-manifold, we'll treat
NQ-maps $T[1]N\map X$ as if they were maps $N\map M$. Likewise, when we
speak about their {\em homotopy}, we mean NQ-maps $T[1](N\times I)\map
X$. If $N$ is a manifold with boundary, an NQ-map $T[1](N\times I)\map
X$ will be called a {\em homotopy rel boundary (r.b.)} if for any
$x\in\partial N$ its restriction to $\{x\}\times T[1]I$ is a constant
map.

A Lie algebroid structure on a vector bundle $A\map M$ is equivalent to
an NQ-structure on the N-manifold $A[1]$. Notice that, according to our
point of view, the groupoid $\Ga$ of \S1 is really a fundamental
groupoid.

NQ-manifolds can have interesting higher homotopies. To make a simple
estimate, let us define the {\em degree} of an N-manifold $X$ (a bit
imprecisely, it the highest degree of a coordinate of $X$). Let $x\in
0\cdot X$ (notice that $0\cdot X$ is an ordinary manifold, as $-1$ acts
trivially there); $x$ does not move under the action of $(\rc,\times)$,
so that the semigroup acts on $T_x X$; the highest weight of this action
is the {\em degree at $x$}. It is locally constant in $0\cdot X$. If
$\deg X=0$ then $X$ is just an ordinary manifold (with trivial
$(\rc,\times)$-action). If $X$ is of degree 1 then it is of the form
$A[1]$ for some vector bundle $A$; hence, an NQ-manifold of degree 1 is
the same as a Lie algebroid.

Here is the estimate:

{\em If $X$ is an NQ-manifold of degree $d$, then $X$ is ``locally a
homotopy $d$-type'', i.e. for any $n>d$, small NQ-maps $T[1]S^n\map X$
can be extended to small NQ-maps $T[1]B^{n+1}\map X$.}

See Lemma 2 in the appendix for a more general (and more appropriate)
formulation.

We can capture the local homotopy of an NQ-manifold using its
fundamental $n$-groupoid, for $n\ge\deg X$: objects are points,
1-morphisms paths connecting the points, 2-morphisms disks connecting
the paths, etc., and finally $n$-morphisms are homotopy classes of
$n$-dimensional balls connecting $n-1$-morphisms (recall that according
to our principle, paths are NQ-maps $T[1]I\map X$ etc.; thus, e.g.,
points are just points in the manifold $0\cdot X$). This $n$-groupoid is
a reasonable generalization of the groupoid corresponding to a Lie
algebroid (for a precise definition of $n$-groupoids via multisimplicial sets, and of fundamental
$n$-groupoids in particular, see \cite{tamsa}). In what follows, we'll
be also interested in various other objects that one can build using
higher homotopies of NQ-manifolds (an example with nice pictures is the
integration of a Lie bialgebroid to a double symplectic groupoid).

\section{Examples: gerbes, loop groups, group-valued moment maps, etc.}

In this section we'll be interested in NQ-manifolds very close to
$T[1]M$, namely in principal $\rc[n]$ bundles over $T[1]M$, $n\ge1$
($\rc[n]$ is $\rc$ or $\rc^{0|1}$, according to the parity of $n$, with
the $(\rc,\times)$-action given by $x\mapsto \la^n x$, and with $Q=0$).
The degree of such an $X$ is $n$ and $0\cdot X=M$.

{\em Lemma: Principal $\rc[n]$-bundles over $T[1]M$ (in the category of
NQ-manifolds) are classified by $H^{n+1}_{DR}(M)$.}

Indeed, in the category of N-manifolds, principal $\rc[n]$-bundles are
trivializable, so let us choose an N-trivialization
$X=T[1]M\times\rc[n]$. The vertical part of $Q$ on $X$ can be identified
with an N-map $T[1]M\map\rc[n+1]$, i.e. with an $n+1$-form $\eta$ on
$M$. Clearly, $Q^2=0$ iff $d\eta=0$. Finally, changing the
trivialization is equivalent to choosing an N-map $T[1]M\map\rc[n]$,
i.e. an $n$-form $\al$; $\eta$ changes to $\eta+d\al$.

Let us also notice that when we choose an N-trivialization of $X$, an
NQ-map $T[1]N\map X$ (where $N$ is a manifold) is the same as a map
$f:N\map M$ and an $n$-form $\om$ on $N$ such that $d\om=f^*\eta$. We
can see easily the $n$-groupoid corresponding to $X$: objects are points
in $M$, 1-morphisms paths in $M$, etc., up to $n-1$-morphisms;
$n$-morphisms  are not just homotopy classes of balls in $M$ rel
boundary (as they would be for $X=T[1]M$), but instead they are points
in a principal bundle over this space of homotopy classes, with the
structure group $\rc/\{\mbox{periods of }\eta\}$ (notice, however, that
$\{\mbox{periods of }\eta\}$ can be dense in $\rc$; in that case the
$n$-groupoid is rather degenerate (recall: beware of large
homotopies)).

Given a principal $U(1)$-bundle over $M$, its Atiyah Lie algebroid gives
us a principal $\rc[1]$-bundle over $T[1]M$. Vice versa, if the periods
of such a $\rc[1]$-bundle are integers, we can integrate it to a
groupoid, whose (any) fibre is a principal $U(1)$ (or $\rc$, if all the
periods vanish) bundle over the universal cover of $M$. There is a
similar relation between $U(1)$-gerbes and $\rc[2]$-bundles, etc.

\subsection{Loop groups}

This was a bit abstract and/or trivial, so let us pass to a nice
concrete (and well studied) example. Let $G$ be a simple compact Lie
group with a chosen invariant inner product. The NQ-group $T[1]G$ has a central extension
$1\map\rc[2]\map\widetilde{T[1]G}\map T[1]G\map1$ that is easily
described in terms of its Lie algebra. Namely, the Lie algebra of
$T[1]G$ is $\fg\ds\fg[1]$ with zero commutator in $\fg[1]$; $Q$ maps
$\fg[1]$ identically to $\fg$, and $\fg$ to $0$. The Lie algebra of
$\widetilde{T[1]G}$ is $\fg\ds\fg[1]\oplus\rc[2]$, where the commutator
in $\fg[1]$ is the inner product; $Q$ acts on the $\fg\ds\fg[1]$-part as
before, and sends $\rc[2]$ to $0$. This defines $\widetilde{T[1]G}$ as
an NQ-group.

Notice that if $M$ is a manifold, the space of NQ-maps
$T[1]M\map\widetilde{T[1]G}$ is a group, and so is the space of
homotopy classes of such maps.

{\em Claim: The group of homotopy classes r.b. of disks in
$\widetilde{T[1]G}$ coincides with the standard central extension
$\widetilde{LG}$ of the loop group.}

Here, as before, a disk in $\widetilde{T[1]G}$ means an NQ-map
$T[1]\bigcirc\map\widetilde{T[1]G}$ (where $\bigcirc$ is a disk).
The claim can be proved by an explicit calculation as follows.
Notice that $G\ss\widetilde{T[1]G}$ (since
$G=0\cdot\widetilde{T[1]G}$). The $\rc[2]$-bundle
$\widetilde{T[1]G}\map T[1]G$ admits a both-sides $G$-invariant
trivialization in the N-category that can be obtained by
exponentiating the inclusion
$\fg\ds\fg[1]\ss\fg\ds\fg[1]\oplus\rc[2]$. A closed 3-form $\eta$
should emerge from this trivialization, and indeed, it is the
3-form coming from the inner product on $G$, $\eta(u,v,w)=\li
u,[v,w]\ri$. Therefore, an NQ-map $T[1]M\map\widetilde{T[1]G}$ is
the same as a map $f:M\map G$ and a 2-form $\om$ s.t.
$d\om=f^*\eta$. One can also check by differentiation that the
product of two NQ-maps is given by $f=f_1f_2$, $\om=\om_1+\om_2+\li
f^*_1\theta_l,f^*_2\theta_r\ri$, where $\theta_{l(r)}$ is the left
(right) Maurer--Cartan form. Now the claim becomes the standard
definition of $\widetilde{LG}$.

\subsection{Moment maps}

The group $\widetilde{T[1]G}$ is also connected with group-valued
moment maps. Namely, let $M$ be a $G$-manifold, so that $T[1]G$
acts on $T[1]M$. $T[1]G$ also acts on $\widetilde{T[1]G}$ by
conjugation. Let us consider a $T[1]G$-equivariant map
$T[1]M\map\widetilde{T[1]G}$. As before, it can be described as a
$G$-equivariant map $f:M\map G$, a 2-form $\om$ s.t. $d\om=f^*\eta$
plus one more condition on $\om$ coming from the
$T[1]G$-equivariance. This is precisely the definition of
$G$-valued moment map of \cite{amm}, except for non-degeneracy of
$\om$. Likewise, if $T[1]M_{1,2}\map\widetilde{T[1]G}$ are two NQ
maps then their product $T[1](M_1\times M_2)\map\widetilde{T[1]G}$
corresponds to the fusion product from the theory of $G$-valued
moment maps.

Now we can describe a simple picture of general moment maps,
following \cite{akm}. Let $F$ be a group with invariant inner
product and $G_{1,2}\ss F$ two Lagrangian subgroups
(=half-dimensional \& isotropic); the case of equivariant moment
maps corresponds to $G_1=G_2$. Notice that
$T[1]G_{1,2}\ss\widetilde{T[1]F}$, as the $G$'s are isotropic. Let
$M$ be a $G_1$-space; consider a $T[1]G_1$-equivariant NQ-map
$T[1]M\map\widetilde{T[1]F}/T[1]G_2$. Again, we could describe it
using a trivialization of the $\rc[2]$-bundle
$\widetilde{T[1]F}/T[1]G_2\map T[1](F/G_2)$; it would become a
$G_1$-equivariant map $f:M\map F/G_2$, a 2-form $\om$, plus some
conditions on $\om$.

Such maps can reasonably be considered as moment maps of general type.
$G$-valued moment maps are included in a simple way: $F=G\times G$, the
inner product is the difference of the inner products on the factors,
$G_1=G_2$ is the diagonal; we have $G\ss F$ (the first factor)
preserving the inner product, and thus also
$\widetilde{T[1]G}\ss\widetilde{T[1]F}$. Now we just identify
$\widetilde{T[1]F}/T[1]G_2$ with $\widetilde{T[1]G}$.

As another example, suppose $G'\ss F$ is another Lagrangian subgroup,
transversal to both $G_{1,2}$. We have $\rc[2]\times
T[1]G'\ss\widetilde{T[1]F}$ ($G'$ is isotropic, so the central extension
becomes trivial there), and we can identify (at least locally)
$\widetilde{T[1]F}/T[1]G_2$ with $\rc[2]\times T[1]G'$. A
$T[1]G_1$-equivariant NQ-map $T[1]M\map\widetilde{T[1]F}/T[1]G_2$ thus
becomes  $G_1$-equivariant map $M\map G'$ (the $G_1$ action on $G'$
comes from the identification of $F/G_2$ with $G'$), a closed 2-form,
and a condition expressing  the equivariance; this is the Lu-Weinstein
moment map. In particular, if $F=T^*G_1$, and $G'=\fg_1^*$, we have the
usual moment map (equivariant if $G_1=G_2$, or twisted by a cocycle if
not).

\subsection{Variational problems}

A reasonably complete understanding of variational problems requires
symplectic geometry (see \S5), but we can make few simple observations
even now.

Let $M$ be a manifold and $N$ an $n$-dim manifold, and let us consider a
Lagrangian for maps $N\map M$ that assigns a density on $N$ to any map
$N \map M$; the action of $f$ is then the integral of this density (for
simplicity, we shall always suppose that $N$ is oriented, and we won't
distinguish between densities and $n$-forms). We can say that the
Lagrangian lifts maps $N \map M$ to  NQ-maps $T[1]N \map
T[1]M\times\rc[n]$.
 As usual, we shall suppose that the density at any
$x\in N$ depends only on the $k$'th jet of the map $N\map M$ at $x$,
where $k$ is some fixed number.

However, as it often happens, a Lagrangian is defined naturally only up
to  a closed $n$-form on $M$. In fact, sometimes it is defined only
locally, with a \v Cech 1-cocycle of closed $n$-forms on $M$ giving the
differences on the overlapping patches. Therefore, we should more
properly say that the Lagrangian lifts maps $f:N\map M$ to NQ-maps $\ti
f:T[1]N\map X$, where $X\map T[1]M$ is a principal $\rc[n]$-bundle. The
action of $f$ is then the homotopy class of $\ti f$ r.b.

For example, in the case of WZW model on a group $G$,
$X=\widetilde{T[1]G}$. The fact that $\widetilde{T[1]G}$ is a
group, plays an  important role again: The solutions of the WZW
model are $g(z,\bar z)=g_1(z)g_2(\bar z)$. The maps
$g_{1,2}:C_{1,2}\map G$ (where $C$'s are (complex) curves) have
unique lifts to NQ-maps $\ti
g_{1,2}:T[1]C_{1,2}\map\widetilde{T[1]G}$ (just because $C$'s are
curves); as it turns out,
 $\ti g=\ti g_1\ti g_2$.

Finally, let us consider symmetries of variational problems.
According to our principle (\S2), the analogue of a 1-parameter
group action $\rc\times M\map M$ is an action $T[1]\rc\times X\map
X$ that is an NQ map. Infinitesimally, such an action is given by
two vector fields $\iota$ and $u$, of degrees $-1$ and $0$
respectively, such that $u=[Q,\iota]$ (so that, after all, we only
need to know $\iota$) and $[\iota,\iota]=0$. Let now, as above, $X$
be a principal $\rc[n]$-bundle over $T[1]M$, and let us consider
only actions of $T[1]I$ preserving the bundle structure (i.e.
commuting with the $\rc[n]$-action). Choosing a (local) NQ
trivialization $X=\rc[n]\times T[1]M$, $\iota$ can be encoded as a
pair $(v,\al)$, where $v$ is a vector field and $\al$ an
$n-1$-form, both on $M$. The condition $[\iota,\iota]=0$ means that
$v\ip\al=0$.

Indeed, pairs $(v,\al)$ (regardless of the condition $v\ip\al=0$)
appear as symmetries of $n$-dim Lagrangians in the usual
formulation. Namely, the pair is a symmetry of $\La$, if ${\cal
L}_v\La+d\al=0$. It leads to a conservation law: if $f:N \map M$ is
extremal, $d(\be_v+f^*\al)=0$, where $\be_v$ is certain $n-1$-form
on $N $ (the momentum density with respect to $v$). Notice that the
conserved momentum $\be_v+f^*\al$ depends on $\al$, not just on
$d\al$. One could expect the pairs $(v,\al)$ to form a Lie algebra,
but the bracket one obtains (the way how $(v_1,\al_1)$ acts on
$(v_2,\al_2)$, namely
$[(v_1,\al_1),(v_2,\al_2)]\equiv([v_1,v_2],{\cal
L}_{v_1}\al_2-v_2\ip d\al_1)$) is not skew symmetric for $n>1$, and
it only gives a Leibniz algebra. This is a small puzzle, but we can
understand it easily now: the bracket is just
$[u_1,\iota_2]=[[Q,\iota_1],\iota_2]$; since
$[[Q,\iota],\iota]=1/2[Q,[\iota,\iota]]$, it may be non-skew if
$\deg X\ge 2$, when $[\iota,\iota]$ can be non-zero. For a complete
description of symmetries we should pass to the differential graded
Lie algebra of $\rc[n]$-invariant vector fields on $X$.

\section{$\Si_n$-manifolds and their homotopy}

A {\em $\Si_n$-manifold} is an NQ-manifold with a $Q$-invariant
symplectic form
 of degree $n$. A {\em $\La$-structure} is a Lagrangian
NQ-submanifold. This section is based on \cite{aksz}.

If $X$ is a $\Si_n$-manifold then $\deg X\le n$. Indeed, if $x\in 0\cdot
X$, any weight $k$ of the $(\rc,\times)$-action on $T_x X$ appears with
the same multiplicity as $n-k$, so that $k\le n$.

A $\Si_0$-manifold is just a symplectic manifold and its
$\La$-structures are Lagrangian submanifolds. A $\Si_1$-manifold is
necessarily of the form $T^*[1]M$ where $M=0\cdot X$; $Q$ has
unique homogeneous Hamiltonian of degree 2 (or of degree $n+1$ for
$\Si_n$-manifolds), which is therefore a bivector field on $M$. In
this way, a $\Si_1$-manifold is the same as a Poisson manifold.
$\La$-structures are just conormal bundles of coisotropic
submanifolds. Similarly, a $\Si_2$-manifold is the same as a
Courant algebroid. If $Y\ss X$ is a $\La$-structure and moreover
$0\cdot Y=0\cdot X$ then $Y$ is the same as a Dirac structure of
the Courant algebroid.\footnote{Courant algebroids were defined by
Z.J. Liu, A. Weinstein and P. Xu \cite{lwx} to serve as Drinfeld
doubles of Lie bialgebroids. Very different (and conceptually
clear) doubles were defined A. Vaintrob (unpublished), using
$\Si_2$-manifolds (with double grading). The connection between the
two approaches was clarified by D. Roytenberg \cite{roy} and A.
Weinstein (unpublished).}

Let now $M$ be a compact oriented $n$-dim manifold, possibly with
boundary, and $X$ a $\Si_n$-manifold, with symplectic form $\om$. The
superspace of all maps $T[1]M\map X$ is symplectic, with symplectic form
denoted $\om_M$, defined as follows: If $\psi:T[1]M\map X$ is any map,
and $u,v\in\Ga(\psi^*TX)$ (i.e. they are infinitesimal deformations of
$\psi$), we let
 $$\om_M(u,v)=\int_{T[1]M}\om(u(y),v(y))dy,$$
 where
$dy$ is the natural Berezin volume form on $T[1]M$ (a function on
$T[1]M$ is a differential form on $M$ and the integral over $T[1]M$ is
the integral of its top-degree part over $M$). Since we chose $\dim
M=n$, we have $\deg\om_M=0$, so that the subspace of all N-maps is
symplectic as well. If we restrict $\om_M$ to the space of NQ-maps, it
is no longer symplectic. However (according to Lemma 3 of Appendix), two
NQ-maps lie on the same null leaf of $\om_M$ iff they are homotopic
r.b.; therefore, the space of homotopy classes r.b. {\em is} symplectic.
Similarly, if we choose a $\La$ structure $Y\ss X$ and consider maps
$T[1]M$ sending $T[1]\partial M$ to $Y$, the space of their homotopy
classes is again symplectic. Finally, if $M=\partial M'$ for some $M'$,
the homotopy classes of the NQ maps extendible to $T[1] M'$ form a
Lagrangian submanifold (again Lemma 3), and in this way we get a
symplectic analogue of  $n+1$-dim topological field theory.
 One should
interpret this with care, since it's based on formal computation of
tangent spaces, cf. the beginning of the appendix.

Let us now pass to examples. The simplest $\Si_2$-manifolds are of the
form $\fg[1]$, where $\fg$ is a Lie algebra with invariant inner product
(any $\Si_2$-manifold $X$ s.t. $0\cdot X$ is a single point is of this
form). If $M$ is a closed oriented surface, the NQ-maps
$T[1]M\map\fg[1]$ are just flat $\fg$-connections on $M$, and the space
of their homotopy classes is the moduli space of these connections, with
its standard symplectic form. If $M$ is a surface with boundary then the
space of homotopy classes r.b. is the space of flat connections modulo
gauge transformations vanishing at the boundary. In particular, if $M$
is a disk, the space is $LG/G$, which is a coadjoint orbit of
$\widetilde{LG}$ (in fact, $\fg[1]$ is a coadjoint orbit of
$\widetilde{T[1]G}$).

Let now $\fh_{1,2}\ss\fg$ be a Manin triple, so that $\fh_{1,2}[1]$ are
mutually transversal $\La$-structures in $\fg[1]$. There is a simple way
of constructing the corresponding double symplectic groupoid (the
symplectic analogue of the quantum group). The groupoid itself is the
space of homotopy classes of maps from this disk to $\fg[1]$:
$$\epsfxsize 3cm \epsfbox{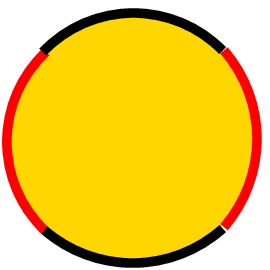}$$
 The black parts of the boundary
are constrained to be mapped to $\fh_1[1]$ and the red parts to
$\fh_2[1]$. (Again, speaking about maps form a disk, we mean NQ
maps from $T[1]$ of the disk, etc.) The two multiplications are
Lagrangian submanifolds of the third Cartesian power of the
groupoid, and they are given by the following picture (this is one
of them---to get the other, just exchange red and black):
 $$\epsfxsize 5cm \epsfbox{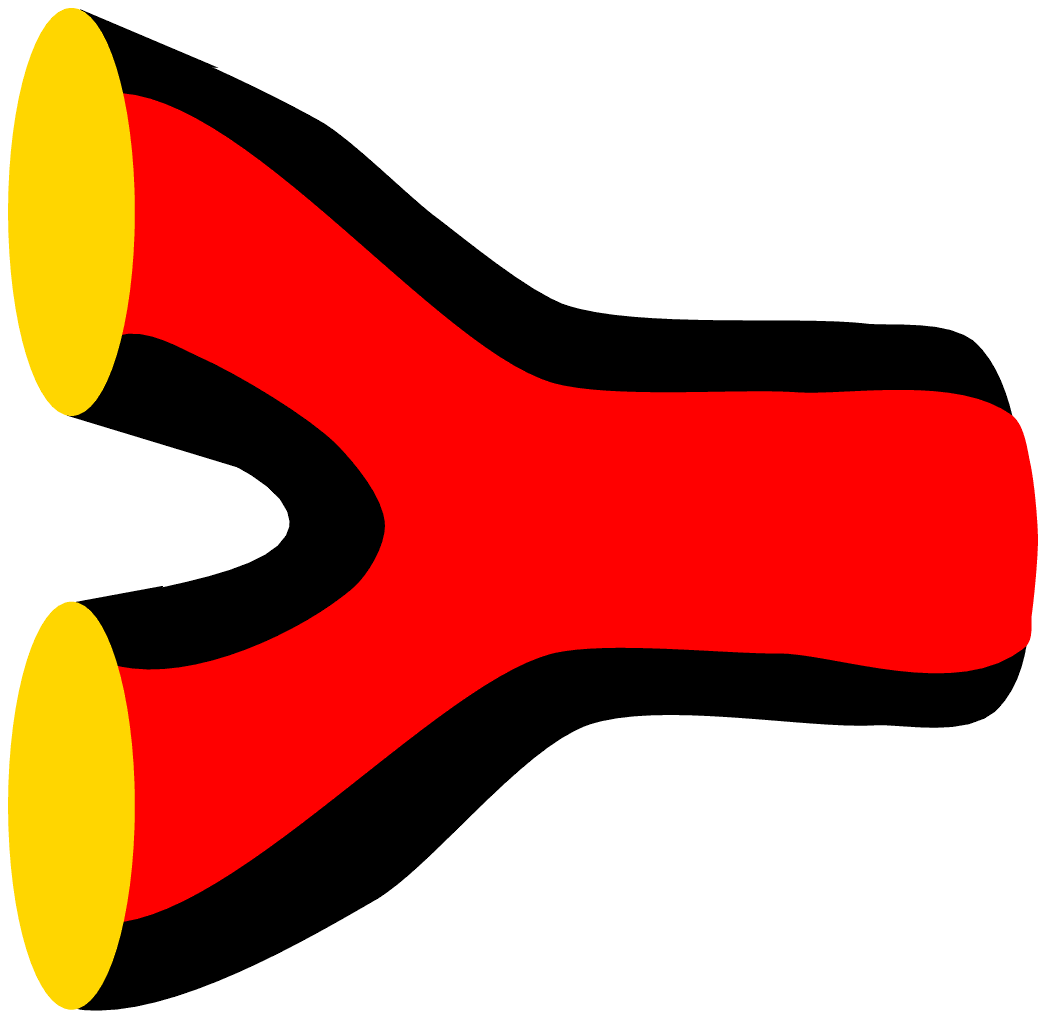}$$
 The
Lagrangian submanifold consists of the homotopy classes of those
maps from three disks that can be extended to a map from the solid
Y on the picture (with the boundary conditions as before). (Other
operations, including the pictures of the Drinfeld double, can be
found in \cite{kwd}.)

The same pictures  can be used for integration of any Lie bialgebroid to
a (local) double symplectic groupoid. We just need to recall the
definition of Lie bialgebroid. It is a $\Si_2$-manifold $X$ with two
$\La$-structures $Y_{1,2}$ such that $0\cdot Y_{1,2}=0\cdot X$ and
$Y_1\cap Y_2=0\cdot X$. The procedure works as before, replacing
$\fg[1]$ with $X$ and $\fh[1]$s with $Y$'s.

This simple construction has several variants. For example, we can turn
one of the red arcs of the circle to pink, and consider one more
$\La$-structure $Z\ss X$ s.t. $0\cdot Z=0\cdot X$; the new boundary
condition is that pink is mapped to $Z$. One of the two multiplications
survives recolouring, so that we still get a symplectic groupoid (though
not double). According to Drinfeld \cite{drinf}, the quadruple
$(X,Y_1,Y_2,Z)$ gives rise to a Poisson homogeneous space; we have just
constructed its symplectic groupoid.

We only dealt with $\Si_2$-manifolds in our examples; many interesting
things may be expected for higher $n$'s.

\section{$\Si_n$-manifolds and variational problems\\(higher-dimensional Hamiltonian mechanics)}

In this section we shall see that $\Si_n$-manifolds play much the same
role in $n$-dim variational problems as Poisson manifolds do in
classical mechanics. In classical mechanics, variational problems lead
to Hamiltonian mechanics in (possibly twisted) cotangent bundles. More
general Poisson manifolds appear most easily by reductions (or sometimes
as degenerate limits). Nevertheless, Hamiltonian mechanics is
interesting even if it doesn't arise in this way. We should mention one
minor problem: 1-dim variational problems lead to Hamiltonian mechanics
only under some invertibility condition (not satisfied e.g. for
reparametrization-invariant Lagrangians). This can be circumvented by
replacing the Hamiltonian---a function on a Poisson manifold $P$---with
a Legendrian submanifold in the space of 1-jets of functions $J^1P$.

One can formulate $n$-dim generalization of Hamiltonian mechanics in the
following way: Let $Y$ be a $\Si_n$-manifold (recall that a
$\Si_1$-manifold is the same as a Poisson manifold). We'll be looking at
NQ-maps $T[1]N\map Y$ (where $N$ is an $n$-dim manifold) satisfying
certain conditions (Hamilton equations). $\La$-structures in $Y$ will
play the role of boundary condition, if we require (parts of)
$T[1]\partial N$ to be mapped there.

Hamilton equations look as follows: let $x\in N$; let us choose a
non-zero element of $\bigwedge^nT_x N$, so that we have a Berezin
integral on $T_x[1]N$. Using this integral, the space of N-maps
$T_x[1]N\map Y$ (denoted $\mbox{\it NMap}(T_x[1]N,Y)$) is a symplectic
manifold (finite-dimensional!). Its symplectic form depends linearly on
the chosen element of $\bigwedge^n T_x N$, i.e. it is a symplectic form
with values in the line $\bigwedge^n T^*_x N$. By definition, the {\em
Hamiltonian$'$} at $x$, denoted $H'_x$, is a Lagrangian submanifold
there. If $T[1]N\map Y$ is an NQ-map, {\em Hamilton equations} require
its restriction to any  $T_x[1]N$ to be in $H'_x$. In the case $n=1$,
$Y=T^*[1]P$ for some Poisson manifold $P$, and $\mbox{\it
NMap}(\rc[1],T^*[1]P)$ is just $T^*P$. As stated above, Hamiltonian is a
Legendrian submanifold in $J^1P$; $H'$ is defined as its projection to
$T^*P$, and we are back in the usual Hamiltonian mechanics.

For general $n$, the {\em Hamiltonian} at $x\in N$, denoted $H_x$, is a
Legendrian submanifold in certain contact manifold $C_x$ (a
generalization of $J^1P$) which is a contactification of $\mbox{\it
NMap}(T_x[1]N,Y)$. $C_x$ is defined as follows: $\mbox{\it
NMap}(T_x[1]N,Y)$ inherits the $(\rc,\times)$ action from $Y$, and
therefore its symplectic form is exact in a canonical way (just plug the
Euler vector field in); we just set $C_x=\mbox{\it
NMap}(T_x[1]N,Y)\times\bigwedge^n T^*_x N$, with the contact structure
given by the $\bigwedge^n T^*_x N$-valued 1-form on $\mbox{\it
NMap}(T_x[1]N,Y)$. Once $H_x$ is chosen, $H'_x$ is defined as its
projection to $\mbox{\it NMap}(T_x[1]N,Y)$. Alternatively, one can
define $C_x$ by first contactifying $Y$ to $Y^C=Y\times\rc[n]$, setting
$C_x=\mbox{\it NMap}(T_x[1]N,Y^C)$ and describing the way in which $C_x$
inherits its contact structure from the contact structure on the
NQ-manifold $Y^C$.\footnote{If $Z$ is a contact N-manifold of degree $n$
and $V$ an $n$-dim vector space then $\mbox{\it NMap}(V[1],Z)$ is a
contact manifold: a vector $v\in T_f(\mbox{\it NMap}(V[1],Z))$ is in the
contact hyperplane at $f$ if $v$, when viewed as a section of $f^*TZ$,
is actually a section of $f^*CZ$, where $CZ\ss TZ$ is the contact
distribution on $Z$. Besides, contact NQ-manifolds are interesting
objects, e.g. those of degree 1 are equivalent to manifolds with Jacobi
structure.}

Let us now describe how this picture arises from variational problems.
We'll be interested in a Lagrangian that to any map $N\map M$ (where $M$
is some manifold) associates a density on $N$; moreover, we'll suppose
that the Lagrangian is first order (i.e. the density depends on the
first jet of the map $N\map M$ only). As we'll see, this situation is
equivalent to the Hamiltonian mechanics described above, when we take
$Y=T^*[n]T[1]M$; Hamiltonian and Lagrangian will become (after some
identifications) the same object. As we discussed in \S3.3, Lagrangian
is natural only up to closed $n$-forms on $M$, and one should use a
principal $\rc[n]$-bundle $X\map T[1]M$ to make it completely natural.
In this picture we get $Y^C=j^1X$ (the space of first jets of sections
of $X\map T[1]M$) and $Y=T^*[n]X//_1\rc[n]$ (the symplectic reduction at
moment 1). One can also go beyond Lagrangians for maps $N\map M$ and
consider an $\rc[n]$-bundle $X\map Z$ over an arbitrary NQ-manifold $Z$
of degree at most $n$ (this is convenient in various kinds of gauge
theories); one still gets $Y^C=j^1X$ and $Y=T^*[n]X//_1\rc[n]$ (if $X$
is trivial then $Y=T^*[n]Z$).

\sloppy

This is how it happens: In the usual terms, a Lagrangian at $x\in
N$ is a $\bigwedge^n T^*_x N$-valued function on the space of
linear maps $T_x N\map TM$, i.e. on the space $\mbox{\it
NMap}(T_x[1]N,T[1]M)$. If we use a principal $\rc[n]$-bundle $X\map
T[1]M$, a Lagrangian at $x$ is a section of the principal
$\bigwedge^n T^*_x N$-bundle
 $\mbox{\it
NMap}(T_x[1]N,X)\map \mbox{\it NMap}(T_x[1]N,T[1]M)$. Finally, to
include multivalued and/or partially defined Lagrangians, we should
define them as Legendrian submanifolds in the space of 1-jets of
sections of this bundle, $j^1\mbox{\it NMap}(T_x[1]N,X)$. We just
notice the natural isomorphism $j^1\mbox{\it
NMap}(T_x[1]N,X)=\mbox{\it NMap}(T_x[1]N,j^1X)$; in this way we get
to Hamiltonian mechanics, with $Y^C=j^1X$ and with Hamiltonian
equal (using the isomorphism) to the Lagrangian.

 \fussy

One question concerning our higher-dimensional Hamiltonian
mechanics certainly remains: {\em what is the meaning of all that?}
We give a rather incomplete answer (incomplete${}={}$not quite
true; please take it with a heap of salt); perhaps this can be
excused by the importance of the question. We pass from $N$ to any
$n$-dim submanifold with boundary (and corners) $N_1\ss N$; the
space of homotopy classes r.b. of NQ-maps $T[1]N_1\map Y$ is
symplectic (as we saw in \S5) and those obeying Hamilton equations
form a Lagrangian submanifold. When we decompose $N_1$ as $N_2\cup
N_3$, these Lagrangian submanifolds compose in the natural way, so
that it's sufficient to know them  for infinitesimal bits of $N$.
The Hamiltonian at $x\in N$ encodes such a Lagrangian submanifold
for an infinitesimal bit around $x$. This picture could be quite
interesting if it survives quantization.

\section{Example: geometry of non-abelian conservation laws in two dimensions}

A conservation law in an $n$-dim variational problem means that we
are given an $n-1$-form which is closed if our configuration is
extremal. For $n=2$ there is a natural non-abelian generalization:
instead of 1-forms, we can consider $\fg$-connections (for some Lie
algebra $\fg$) that are flat on extremals. One can consider
non-abelian generalizations for higher $n$'s as well: a flat
connection is an NQ-map to $\fg[1]$, and in principle we can
replace $\fg[1]$ with any NQ-manifold. However, the geometry of
these general conservation laws is unclear to me; we shall restrict
ourselves to $n=2$ and to $\fg[1]$.

Thus, we shall consider maps from a surface $N$ to a manifold $M$;
any such map should produce a $\fg$-connection on $N$, flat if the
map is extremal. More naturally, we'll be given a principal
$G$-bundle $P\map M$; any map $f:N\map M$ should produce a
connection on $f^*P$, flat if $f$ is extremal. We will call this a
{\em  $P$-conservation law}. A flat connection on $f^*P$ is the
same as a lift of $f$ to an NQ-map $T[1]N\map(T[1]P)/G$.

\sloppy

The idea is as follows: The variational problem will be
reformulated  as 2dim Hamiltonian mechanics in a $\Si_2$-manifold
$Y_M$. We shall then define a coisotropic NQ-submanifold
$Y_{M|0}\ss Y_M$ and an NQ-map $Y_{M|0}\map(T[1]P)/G$. If the
Hamiltonian is such that the Hamilton equations force maps
$T[1]N\map Y_M$ to be maps $T[1]N\map Y_{M|0}$ then we have a
$P$-conservation law. Since $Y_{M|0}$ is coisotropic, this
condition is easily understood. Recall that for any $x\in N$,
$H'_x$ is a Lagrangian submanifold in $\mbox{\it
NMap}(T_x[1]N,Y_M)$; we want it to be a submanifold of  $\mbox{\it
NMap}(T_x[1]N,Y_{M|0})$. Since $\mbox{\it NMap}(T_x[1]N,Y_{M|0})$
is a coisotropic submanifold of $\mbox{\it NMap}(T_x[1]N,Y_M)$, it
just means that $H'_x$ is a maximally isotropic submanifold of the
presymplectic manifold $\mbox{\it NMap}(T_x[1]N,Y_{M|0})$; in
particular, it is woven of the null leaves of the presymplectic
form.  In the favorable case, when the space of null leaves is a
supermanifold (hence a $\Si_2$-manifold), it reduces to a
Hamiltonian on this space of leaves.

\fussy

$Y_{M|0}\ss Y_M$ and $Y_{M|0}\map (T[1]P)/G$ will be constructed
via a group action, and the condition on Hamiltonian will be
reformulated as its invariance, i.e. we get a non-Abelian version
of Noether theorem. Namely, consider a Lie algebra $\ff$ with
invariant inner product, containing $\fg$ as a Lagrangian
subalgebra (i.e. $\fg\ss\ff$ is a Manin pair). We shall suppose
that $F$ acts on $P$, extending the action of $G$ (the action of
$F$ needn't preserve the bundle structure). In addition, we shall
also need a principal $\rc[2]$-bundle $X_P\map T[1]P$ (in
NQ-category) and an NQ-action of $\widetilde{T[1]F}$ on $X_P$,
covering the action of $T[1]F$  on $T[1]P$ and extending the action
of $\rc[2]$ on $X_P$. We let $X_M=X_P/T[1]G$ (thus $X_M\map T[1]M$
is a principal $\rc[2]$-bundle). As in \S5, we define $Y_{P,M}$ as
$T^*[2]X_{P,M}//_1\rc[2]$. Clearly $Y_M=Y_P//T[1]G$.

The group $\widetilde{T[1]F}$ acts  on $Y_P$; as
$\rc[2]\ss\widetilde{T[1]F}$ acts trivially, the action factors
through a (non-Hamiltonian) $T[1]F$-action. The moment map $\mu$ of
the $\widetilde{T[1]F}$-action takes values in
$T^*[2](\widetilde{T[1]F})/\widetilde{T[1]F}=\rc\oplus\ff^*[1]\oplus\ff^*[2]$;
the $\rc$-component of $\mu$ is identically 1. Let $\mu_0=1 \oplus
0 \oplus 0 \in \rc\oplus\ff^*[1]\oplus\ff^*[2]$; $\mu_0$ is clearly
a fixed point of $Q$ and it is also $(\rc,\times)$-invariant, so
$Y_{P|\mu_0}=\{x\in Y_P, \mu(x)=\mu_0\}\ss Y_P$ is a NQ-manifold
and so is the symplectic reduction
$Y_P//_{\mu_0}\widetilde{T[1]F}=Y_{P|\mu_0}/F$, if it is a
supermanifold ($F$ is the subgroup of $T[1]F$ fixing $\mu_0$).
Notice that $Y_{P|\mu_0}/G\ss Y_M$ is a coisotropic NQ-submanifold
(as $Y_M=Y_P//T[1]G$). We let $Y_{M|0}=Y_{P|\mu_0}/G$. If
$Y_P//_{\mu_0}\widetilde{T[1]F}$ is a supermanifold, it is the
space of null leaves of $Y_{M|0}$. The NQ-map $Y_{M|0}\map
(T[1]P)/G$ needed for $P$-conservation law comes from the
projection $Y_P\map T[1]P$.

To conclude, we should explain how our condition on the Hamiltonian
can be expressed as its symmetry. We define an auxiliary symplectic
N-submanifold $Y_{P|>}$ of $Y_P$ as the subspace where the
$\ff^*[1]$-component of $\mu$ vanishes. The group $F$ acts on
$Y_{P|>}$ and its moment map vanishes on $Y_{P|0}$. Let $H^P_x$ be
a Legendrian submanifold in $NMap(T_x[1]N,Y_{P|>}^C)$; it can be
projected to a Legendrian submanifold in $NMap(T_x[1]N,Y_{M}^C)$,
i.e. to a Hamiltonian in $Y_M$, iff it is $G$-invariant, and the
Hamiltonian will satisfy our condition iff $H^P_x$ is
$F$-invariant.

\sloppy

 Finally, there is one more interesting phenomenon
concerning $P$-conservation laws. Suppose there is another
Lagrangian subgroup $\bar{G}\ss F$ such that $P\map P/\bar{G}$ is a
principal $\bar{G}$-bundle. Then we can define another variational
problem with target $\bar{M}=P/\bar{G}$, just by projecting $H^P_x$
to a Legendrian submanifold of $NMap(T_x[1]N,Y_{\bar{M}}^C)$.
Moreover, one can transfer solutions of Hamilton equations between
$Y_M$ and $Y_{\bar{M}}$: Recall that our Hamilton equation
constrain NQ-maps $T[1]N\map Y_M$ to be NQ-maps $T[1]N\map
Y_{P|\mu_0}/G$. The latter maps can be lifted to NQ-maps $T[1]N\map
Y_{P|\mu_0}$ (provided $\pi_1(N)=0$, otherwise the lift may be
multivalued) and these lifts are unique up to $G$-action, i.e. the
projection $Y_{P|\mu_0}\map Y_{P|\mu_0}/G$ is a ``normal covering
with group $G$" (this is because $G$ is ``discrete" in the world of
NQ-manifolds, i.e. for any manifold $K$, any NQ-map $T[1]K\map G$
is constant). Then we project this lift to an NQ-map $T[1]N\map
Y_{{\bar{M}}|0}$. This procedure takes solutions to solutions and
in this way the two variational problems become equivalent. This is
the {\em Poisson-Lie T-duality} of \cite{ks}.

\fussy

\section*{Appendix: Infinitesimal deformations and vector bundles}

Let $X$ and $Y$ be NQ-manifolds; recall that by a homotopy
connecting two NQ-maps $Y\map X$ we mean an NQ-map $Y\times
T[1]I\map X$. We'll be interesting in the NQ-maps $Y\map X$ close
to a given NQ-map $\psi:Y\map X$, and also in the question which of
them are homotopic to each other by small homotopies. We will
simply linearize the problem (considering infinitesimally close
maps) and compute formally the tangent space to the space of
homotopy classes of NQ-maps $Y\map X$. It is by no means true that
this space is always a manifold and we don't give any criteria for
the formal tangent space to be actual tangent space. This is a
serious flaw in this paper.

An infinitesimal deformation of $\psi:Y\map X$ is a section of
$\psi^*TX$; the infinitesimal part of a homotopy starting at $\psi$
is a linear homotopy connecting a section of $\psi^*TX$ with 0. By
an {\em NQ vector bundle}  $E\map Y$ we mean a vector bundle that
is both $(\rc,\times)$- and $Q$-equivariant ($\psi^*TX$ was an
example). Notice that the vector superspace of its sections is a
cochain complex (a section $s$ is of degree $d$ if $s(\la\cdot
x)=\la^d(\la\cdot s(x))$, $\la\in\rc$, where $\cdot$ denotes the
$(\rc,\times)$-action; $d$ may be negative, but it's bounded from
below, e.g. $d\ge-\deg X$ for $E=\psi^*TX$).  A section is an N-map
if it's of degree 0, it is an NQ-map if it's also closed, and it is
connected with the zero section by a linear homotopy if it's exact.
Notice however that {\em any} homotopy in the space of sections can
be made to a linear homotopy. Hence, $H^0(\Ga(E))$ is the space of
homotopy classes of NQ sections, whether the homotopies are
required to be linear or not.

\sloppy

Negative cohomologies have a homotopical meaning as well; namely,
$H^{-n}(\Ga(E))$ is ``$\pi_n$ of the space of sections''. More
precisely, since $H_{DR}(B^n,S^{n-1})=\rc[n]$, we have:

\fussy

{\em Lemma 1: Let $E\map Y\times T[1]B^n$ be an NQ vector bundle, and
let $E_0$ be its restriction to $Y$ (where $Y$ is embedded into $Y\times
T[1]B^n$ as (say) $Y\times\{\mbox{\rm centre of the ball}\}$). Let
$\Ga_0(E)$ denote the sections of $E$ vanishing at $Y\times
T[1]S^{n-1}$. Then $H(\Ga_0(E))=H(\Ga(E_0))[n]$.}

This lemma can also be generalized to a form of Thom isomorphism
(replacing $Y\times T[1]B^n$ by a bundle).

Now, let $\psi:T[1]B^n\map X$ be an NQ-map. Setting $Y=\{\mbox{\it
pt}\}$, $E=\psi^*TX$, the previous lemma gives us a linearized version
of the following:

{\em Lemma 2: Let $\psi:T[1]B^n\map X$ be an NQ-map. If $n>\deg X$
then any NQ-map close to $\psi$ and coinciding with $\psi$ on
$T[1]S^{n-1}$ is homotopic to $\psi$ r.b. by a small homotopy, i.e.
$X$ is ``locally a homotopy $\deg X$-type''.}

It is possible to deduce this claim from its linearized version (a nice
explanation would be nice).

Let us now pass to symplectic vector bundles. A {\em $\Si_n$ vector
bundle} is an NQ vector bundle with symplectic fibres, with the field of
symplectic forms of degree $n$, and with $Q$ compatible with the
symplectic structure.

We'll be mostly interested in $\Si_n$  bundles $E\map T[1]M$, where $M$
is an oriented compact manifold, possibly with boundary. Notice that
$\Ga(E)$ is a symplectic vector superspace; the symplectic form (denoted
$\om_M$) is given by
 $$\om_M(u,v)=\int_{T[1]M}\om(u(y),v(y))dy,$$
where $dy$ is the standard Berezin volume form on $T[1]M$. The degree of
$\om_M$ is $n-\dim M$.

Let $E'$ be the restriction of $E$ to $T[1]\partial M$ and $u'$, $v'$
the restrictions of $u$ and $v$. Notice (by Stokes theorem) that
$$\om_{\partial M}(u',v')=\pm\om_M(Qu,v)\pm\om_M(u,Qv).$$
 As a
consequence,

{\em Lemma 3: $Z\Ga(E)=(B\Ga_0(E))^\bot$ (where $\Ga_0(E)$ denotes the
subcomplex of sections vanishing at $T[1]\partial M$), so that
$Z\Ga(E)/B\Ga_0(E)$ is symplectic. In particular, if $M$ is closed,
$H(\Ga(E))$ is symplectic (a version of Poincar\'e duality). The image
of $H(E)\map H(E')$ is a Lagrangian subspace.}

Notice that if $\dim M=n$, $\deg\om_M=0$, so that
$(Z\Ga(E)/B\Ga_0(E))^0$, the space of homotopy classes r.b. of
NQ-sections of $E$, is symplectic. Finally notice that if $M$ is not
closed, $H(E)$ is not symplectic. But there is a simple way how to make
$Q$ selfadjoint. We choose a Lagrangian NQ subbundle $L\ss E'$, and let
$\Ga_L(E)$ be the sections taking values in $L$ over $T[1]\partial M$.
Then $H(\Ga_L(E))$ is symplectic.

\newpage
\unitlength 1cm
\begin{picture}(10,19)
\put(1,0){\begin{sideways}
\begin{tabular}{|c|c|c|c|c|c|} \hline
$n$ & &{\em substructures}&{\em integration${}={}$homotopy}&{\em
quantization}&{\em variational problems}\\ \hline\hline
0&symplectic&Lagrangian&---&vector spaces&---\\
&manifolds&submanifolds&&&\\ \hline
1&Poisson&coisotropic&symplectic groupoids&associative algebras,&particle\\
&manifolds&submanifolds&&abelian categories&mechanics\\ \hline
&Courant&(generalized)&symplectic 2-groupoids,&abelian 2-categories,&2-dim variational\\
2&algebroids&Dirac&double symplectic groupoids, &quantum groups,&problems\\
&&structures&3-dim symplectic TFTs, etc.&3-dim TFTs&\\ \hline
$n$&$\Si_n$-manifolds&$\La$-structures&symplectic $n$-groupoids,&abelian $n$-categories,&$n$-dim variational\\
&&&$n+1$-dim symplectic TFTs, etc.&$n+1$-dim TFTs&problems\\ \hline
\end{tabular}
\end{sideways}}
\put(7.5,9){\begin{sideways}\bf Table 1\end{sideways}}
\put(8,7.2){\begin{sideways}\it Whatever this table is
about\end{sideways}}
\end{picture}

\end{document}